\documentstyle[12pt,leqno]{article}

\newtheorem{thm}{Theorem}[section]
\newtheorem{Theorem}{Theorem}[section]

\newcommand{\sthm}{\begin{Theorem}}         %start theorem
\newcommand{\ethm}{\end{Theorem}}           %end theorem
\newtheorem{Corollary}[Theorem]{Corollary}
\newcommand{\scor}{\begin{Corollary}}       %start corollary
\newcommand{\ecor}{\end{Corollary}}         %end corollary
\newcommand{\pf}{ \par \vspace{1ex} \noindent {\sc Proof.} \hspace{2mm}}

\begin{document}
\title{The explicit formula of flat Lagrangian $H$-umbilical submanifolds in quaternion Euclidean spaces}
\author{Yun Myung Oh and Joon Hyuk kang\\
ohy@math.msu.edu, kang@andrews.edu}
\date{ }
\maketitle
\begin{abstract}
In \cite{oh}, there exist nonflat Lagrangian $H$-umbilical
submanifolds in $\bf H^{n}$: Lagrangian pseudo-sphere and a
quaternion extensor of the unit hypersphere of $\bf E^{n}$. In
this paper, using the idea of twisted product, we investigate the
flat Lagrangian $H$-umbilical submanifolds in quaternion Euclidean
space $\bf H^{n}$.
\end{abstract}
\section{Introduction}
We begin with the following results from \cite{c97}. B.Y. Chen
introduced the notion of Lagrangian $H$-umbilical submanifold in
$\bf C^{n}$ and classified the Lagrangian $H$-umbilical sumaifolds
$M$ in $\bf C^{n}$:flat, Lagrangian pseudo-sphere, or complex
extensor of the unit hypersphere of $\bf E^{n}$. Moreover, he also
obtained the following results for the flat case. Let $M$ be a
simply-connected open portion of the twisted product manifold
$_{f}\bf R \times E^{n-1}$ with the metric
$$g=f^{2}dx_{1}^{2}+\sum_{j=2}^{n}dx_{j}^{2},$$ where
$f=\beta(x_{1})+\sum_{j=2}^{n}\alpha_{j}(x_{1})x_{j}$ for some
real valued functions $\beta$ and $\alpha_{2},...,\alpha_{n}$.
There exists a $unique$ Lagrangian isometric immersion
$L_{f}:M\rightarrow C^{n}$ with the second fundamental form
$$h(e_{1},e_{1})=\lambda Je_{1},
h(e_{1},e_{j})=h(e_{j},e_{k})=0, 2\leq j,k\leq n,$$ where
$\lambda=\frac{1}{f},e_{1}=\lambda \frac{\partial} {\partial
x_{1}},e_{2}=\frac{\partial}{\partial x_{2}},
\cdots,e_{n}=\frac{\partial}{\partial x_{n}}$. Furthermore, if
$n\geq 3$, and $L:M \rightarrow C^{n}$ is a Lagrangian
$H$-umbilical isometric immersion of a flat manifold into $C^{n}$
without totally geodesic points, then $M$ is an open portion of a
twisted product manifold $_{f}\bf R \times E^{n-1}$ described as
above. Later, in \cite{c}, B.Y. Chen obtained the explicit
description of this isometric immersion using the idea of
a special Legendre curve.\\
Based on the facts above, we can impose the following questions.\\
\underline{Question 1}. Are there any flat Lagrangian $H$-umbilical
submanifolds in quaternion Euclidean spaces? \\
I could answer this question in theorem 4.1 \cite{oh} and get more
details in section 3. In fact, besides this flat submanifolds,
there exists a Lagrangian pseudo-sphere in $\bf C^{n}$ and a
quaternionic extensor of the unit hypersphere of $\bf E^{n}$.\\
\underline{Question 2}. If exists, what is the explicit description
of this isometric immersion?
\\[0.1in]
We remark here that we follow the notations and definitions
given in \cite{oh}.\\
\section{Preliminaries}
We have the following existence and uniqueness theorems.\\
\begin{thm}\label{Theorem 2.1} Let $M^{n}$ be a simply connected
Riemannian $n$-manifold and $\sigma_{i}(i=1,2,3)$ be  $TM$-valued
symmetric bilinear forms on $M$ such that\\
(a) $<\sigma_{i}(X,Y),Z>$ is totally symmetric for $i=1,2,3$\\
(b) $(\nabla_{X}\sigma_{i})(Y,Z)-\sigma_{j}(X,\sigma_{k}(Y,Z))+
\sigma_{k}(X,\sigma_{j}(Y,Z))$ is totally symmetric, where
$(\nabla_{X}\sigma_{i})(Y,Z)=\nabla_{X} \sigma_{i}(Y,Z)-
\sigma_{i}(\nabla_{X}Y,Z)-\sigma_{i}(Y,\nabla_{X}Z)$ and
$(i,j,k)= (1,2,3), (2,3,1),\mbox{ or } (3,1,2)$.\\
(c) $R(X,Y)Z=\sum_{i=1}^{3}\{\sigma_{i}(\sigma_{i}(Y,Z),X)-
\sigma_{i}(\sigma_{i}(X,Z),Y)\}$\\
Then there exists a Lagrangian isometric immersion
$x:M^{n} \rightarrow \bf {H}^{n}$ whose second fundamental
form $h(X,Y)=I\sigma_{1}(X,Y)+J\sigma_{2}(X,Y)+K\sigma_{3}(X,Y)$.
\end{thm}
\begin{thm}\label{Theorem 2.2} Let $L_{1},L_{2}:M^{n} \rightarrow \bf H^{n}$
be two Lagrangian isometric immersion of a Riemannian $n$-manifold
with the second fundamental forms $h^{1}$ and $h^{2}$,
respectively. If, for $i=1,2,3$
$$<h^{1}(X,Y),\pi_{i}L_{1^{\star}}Z>=<h^{2}(X,Y),\pi_{i}L_{2^{\star}}Z>$$
for all tangent vector fields $X,Y,Z$ on $M^{n}$, and $\pi_{i}=I, J$ or $K$,
then there exists an isometry $\phi$ of $\bf H^{n}$ such that $L_{1}=L_{2} \circ \phi$.
\end{thm}
Now, here is the sketch of the proof of theorem 2.1.\\
\underline{Proof of theorem 2.1} We define a bundle $NM$ over $M$ by
$NM=TM\oplus TM\oplus TM$, $\pi_{i}:TM\rightarrow NM$, where
$\pi_{1}(X)=(X,0,0),\pi_{2}=(0,X,0),\pi_{3}=(0,0,X)$. We also
define a connection on $NM$ by
\[\begin{array}{ccc}
\nabla^{\perp}_{X}(\pi_{1}Y_{1}+\pi_{2}Y_{2}+\pi_{3}Y_{3}) & =
& \pi_{1}\nabla_{X}Y_{1}+\pi_{3}\sigma_{2}(X,Y_{1})-\pi_{2}\sigma_{3}(X,Y_{1})\\
& + & \pi_{2}\nabla_{X}Y_{2}-\pi_{3}\sigma_{1}(X,Y_{2})+\pi_{1}\sigma_{3}(X,Y_{2})\\
& + & \pi_{3}\nabla_{X}Y_{3}+\pi_{2}\sigma_{1}(X,Y_{3})-\pi_{1}\sigma_{2}(X,Y_{3}).
\end{array}\]
Then we can define the second fundamental form $h:TM\times TM
\rightarrow NM$ by
$h(X,Y)=\pi_{1}\sigma_{1}(X,Y)+\pi_{2}\sigma_{2}(X,Y)+\pi_{3}\sigma_{3}(X,Y)$.
Its corresponding Weingarten maps are given by
$A_{\pi_{i}X}Y=\sigma_{i}(X,Y),i=1,2,3$. Then the straightforward
long calculations show that this setting satisfies the Gauss,
Codazzi and Ricci equations. Applying the Existence theorem, there
exists an isometric immersion $x:M^{n}\rightarrow \bf E^{4n}$ with
the normal bundle $NM$, second fundamental form $h$, normal
connection $\nabla^{\perp}$ and Weingarten operator $A$.  Let's
define three endomorphisms $I,J$ and $K$ on $\bf E^{4n}=TM+NM$ as
below:\\[0.1in]
 \( \begin{array}{cccc}
I|_{TM}=\pi_{1}TM & I|_{\pi_{1}TM}=-TM & I|_{\pi_{2}TM}=\pi_{3}TM & I|_{\pi_{3}TM}=-\pi_{2}TM \\
J|_{TM}=\pi_{2}TM & J|_{\pi_{2}TM}=-TM & J|_{\pi_{3}TM}=\pi_{1}TM & J|_{\pi_{1}TM}=-\pi_{3}TM \\
K|_{TM}=\pi_{3}TM & K|_{\pi_{3}TM}=-TM & K|_{\pi_{1}TM}=\pi_{2}TM
& K|_{\pi_{2}TM}=-\pi_{1}TM
\end{array} \)
\\[0.1in]
It is easy to check that these three endomorphisms are almost
complex structures satisfying:
$$I^{2}=J^{2}=K^{2}=-1,IJ=K,JI=-K,JK=I,KJ=-I,KI=J,IK=-J$$
Using these defintions, we can easily verify that the second
fundamental form $h$ is now given by
$$h(X,Y)=I\sigma_{1}(X,Y)+J\sigma_{2}(X,Y)+K\sigma_{3}(X,Y)$$
Finally, we must show that $I,J$ and $K$ are parallel. For $X,Y$
tangent vector fields to $M$, we get
\\[0.1in]
 \( \begin{array}{ccl}
(\tilde{\nabla}_{X}I)Y & = &
-A_{IY}X+\nabla^{\perp}_{X}(IY)-I\nabla_{X}Y-Ih(X,Y) \\
& = &
-\sigma_{1}(X,Y)+\pi_{1}\nabla_{X}Y-\pi_{2}\sigma_{3}(X,Y)+\pi_{3}\sigma_{2}(X,Y) \\
&  & -\pi_{1}\nabla_{X}Y-Ih(X,Y) = 0. \end{array} \)
\\[0.1in]
Similarly, we can also show that $(\tilde{\nabla}_{X}J)Y=0,
(\tilde{\nabla}_{X}K)Y=0$ and
$(\tilde{\nabla}_{X}\varphi)(\pi_{i}Y)=0$, where $\varphi=I,J$ or
$K$ and $i=1,2,3$.  Therefore, we can identify ${\bf E^{4n}}, I,
J$ and $K$ with $\bf H^{n}$ and we can easily see that this
isometric immersion $x$ is Lagrangian.
\\[0.1in]
Now, we recall a definition of twisted product \cite{c97}. Let $N_{1},N_{2}$ be
two Riemannian manifolds with Riemannian metrics $g_{1},g_{2}$,
respectively and $f$ a positive function on $N_{1}\times N_{2}$.
Then the metric $g=f^{2}g_{1}+g_{2}$ is called a twisted product
metric on $N_{1}\times N_{2}$. The manifold $N_{1}\times N_{2}$
with the twisted product metric $g=f^{2}g_{1}+g_{2}$ is called a
twisted product manifold, which is denoted by $ _{f}N_{1}\times N_{2}$.
The function $f$ is called the twisting function of the twisted product manifold.\\
\section{Main results}
In order to characterize the flat Lagrangian submaifold into
quaternion Euclidean spaces, we need the quaternion version of
special Legendre curve in $S^{n-1}\in {\bf C^{n}}$ introduced
by B.Y. Chen in his paper \cite{c}.\\
Let $z:I\rightarrow S^{4n-1}\subset {\bf H^{n}}$ be a unit speed
curve in the unit hypersphere centered at the origin in ${\bf
H^{n}}$ satisfying the following condition: $<z'(s),iz(s)>=
\\
<z'(s),jz(s)>=<z'(s),kz(s)>=0$ identically. Hence
$z(s),iz(s),jz(s),kz(s),\\
z'(s),iz'(s),jz'(s),kz'(s)$ are orthonormal vector fields defined
along the curve. Thus, there exists normal vector fields
$P_{3},P_{4},...,P_{n}$ such that $z(s),iz(s),jz(s),
\\
kz(s),z'(s),iz'(s),jz'(s),kz'(s),P_{3},iP_{3},
jP_{3},kP_{3},...,P_{n},iP_{n},jP_{n},kP_{n}$ form an orthonormal
frame field along the Legendre curve. Using these orthonormal
vector fields, $z''(s)$ can be written as
\begin{eqnarray}  \label{eq:1}
z''(s) & = & i\alpha(s)z'(s)+j\beta(s)z'(s)+k\gamma(s)z'(s)-z(s)-
\sum_{l=3}^{n}a_{l}(s)P_{l}(s) \\
& + & \sum_{l=3}^{n}b_{l}(s)iP_{l}(s)+
\sum_{l=3}^{n}c_{l}(s)jP_{l}(s)+\sum_{l=3}^{n}d_{l}(s)kP_{l}(s),\nonumber
\end{eqnarray}
where $\alpha,\beta,\gamma,a_{l},b_{l},c_{l}$ and $d_{l}$ are all
real valued functions. The Legendre curve $z=z(s)$ is called a
special Legendre curve in $S^{4n-1}\subset {\bf H^{n}}$ if the
expression (\ref{eq:1}) is simplified to
$$z''(s)=i\alpha(s)z'(s)+j\beta(s)z'(s)+k\gamma(s)z'(s)-z(s)-
\sum_{l=3}^{n}a_{l}(s)P_{l}(s)$$ for some parallel normal vector
fields $P_{3},P_{4},...,P_{n}$ along the curve.\\
We note here that any Legendre curve in $S^{7}\subset {\bf H^{2}}$
is special.
\begin{thm}\label{theorem:3.1}(a) Let $M^{n}$ be a simply connected open
portion of the twisted product manifold $_{f}\bf R\times \bf E^{n-1}$
with twisted product metric
$$g=f^{2}dx_{1}^{2}+\sum^{n}_{j=2}dx_{j}^{2}$$
where $f^{2}=f_{1}^{2}+f_{2}^{2}+f_{3}^{2}$ for three arbitrary functions
on $M$ such that $\frac{f}{f_{1}},\frac{f}{f_{2}},\frac{f}{f_{3}}$ are
functions of only $x_{1}$ and also $$f(x_{1},...,x_{n})=\beta(x_{1})+
\sum_{j=2}^{n}\alpha_{j}(x_{1})x_{j}$$ for some functions $\beta,\alpha_{1},...,\alpha_{n}$ of $x_{1}$. \\
Then, up to rigid motions of $\bf H^{n}$, there is a unique Lagrangian
isometric immersion $L_{f}:M \rightarrow \bf H^{n}$ without totally geodesic
points whose second fundamental form satisfies\\
$$h(e_{1},e_{1})=\lambda_{1}Ie_{1}+\lambda_{2}Je_{1}+\lambda_{3}Ke_{1},
h(e_{1},e_{j})=h(e_{j},e_{k})=0\;\;\; j,k\geq 2,$$
where $e_{1}=\displaystyle{\frac{1}{f} \frac{\partial}{\partial x_{1}}},
e_{i}=\frac{\partial}{\partial x_{i}}, i\geq 2,
\lambda_{j}=\frac{1}{f_{j}},j=1,2,3$\\[0.1in]
(b)  Suppose $L:M^{n}\rightarrow \bf H^{n} (n\geq 3)$ is a Lagrangian
$H$-umbilical isometric immersion of a flat manifold into $\bf H^{n}$
without totally geodesic points. Then $M$ is an open portion of a twisted
product $_{f}\bf R \times  \bf E^{n-1}$ with twisted product metric
$g=f^{2}dx_{1}^{2}+dx_{2}^{2}+ \cdots + +dx_{n}^{2}$ and twisted
product fuction given in statement (a). Up to rigid motions of
$\bf H^{n}$, $L$ is the uniquely given by the $L_{f}$ above in (a).
However,\\
\hspace*{.5in}(i) If $\alpha_{2}=\cdots=\alpha_{n}=0$, i.e. $f$ is
a function of $x_{1}$ only, then $L$ is given by
$L(x_{1},...,x_{n})=D(x_{1})+\sum_{j=2}^{n}c_{j}x_{j}$ which is a
Lagrangian cylinder over a curve $D(x_{1})$
whose rulings are $(n-1)$ planes parallel to $x_{2}\cdots x_{n}$-planes in $\bf H^{n}$.\\
\hspace*{.5in}(ii) Otherwise, by doing some change of
variables $t=\int^{x_{1}}_{0} \alpha_{2}(x)dx, u_{2}=
x_{2},\cdots,u_{n}=x_{n}$, $L$ is given by $$L(t,u_{2},
\cdots,u_{n})=u_{2}z(t)+\sum_{j=3}^{n}u_{j}P_{j}(t)+\int^{t}b(t)z'(t)dt$$
for some $\bf H^{n}$ valued functions $P_{3},...,P_{n}$ of $t$,
where the twisted product metric $g=\tilde{f}^{2}dt^{2}+du_{2}^{2}+\cdots+du_{n}^{2},$
and twisted product function $\tilde{f}(t,u_{2},...,u_{n})=b(t)+u_{2}+
\sum^{n}_{j=3}a_{j}(t)u_{j}$. Here, $z=z(t)$ is a special Legendre curve in $S^{4n-1}$.
\end{thm}
\pf (a) Define three symmetric bilinear forms
$\sigma_{1},\sigma_{2}$ and $\sigma_{3}$ on $M^{n}$ by
$\sigma_{1}(\frac{\partial}{\partial
x_{1}},\frac{\partial}{\partial
x_{1}})=\frac{f}{f_{1}}\frac{\partial}{\partial x_{1}},
\sigma_{2}(\frac{\partial}{\partial
x_{1}},\frac{\partial}{\partial
x_{1}})=\frac{f}{f_{2}}\frac{\partial}{\partial x_{1}},
\sigma_{3}(\frac{\partial}{\partial
x_{1}},\frac{\partial}{\partial
x_{1}})=\frac{f}{f_{3}}\frac{\partial}{\partial x_{1}}$ and all
other are zero. Then $<\sigma_{i}(X,Y),Z>$ is totally symmetric in
$X,Y$, and $Z$. Using the twisted product metric given, we have
$$\left.\begin{array}{l}
\nabla_{\frac{\partial}{\partial x_{1}}}\frac{\partial}{\partial
x_{1}}  = \frac{(f)_{1}}{f}\frac{\partial}{\partial
x_{1}}-f\sum_{2}^{n}(f)_{k}\frac{\partial}{\partial x_{k}}, \\
\nabla_{\frac{\partial}{\partial x_{1}}}\frac{\partial}{\partial
x_{i}}  = \frac{(f)_{i}}{f}\frac{\partial}{\partial x_{1}},
\nabla_{\frac{\partial}{\partial x_{j}}}\frac{\partial}{\partial
x_{k}}=0,i,j,k=2,...,n,
\end{array}\right.$$
where $(f)_{j}=\frac{\partial f}{\partial x_{j}}$ for $j=1,...,n$.
We note here that
$(f)_{1}=\beta'(x_{1})+\sum_{j=2}^{n}\alpha'_{j}(x_{1})x_{j}$, and
$(f)_{j}=\alpha_{j}(x_{1})$ for $j=2,...,n$. The long
straightforward computations show that all three conditions of
theorem 2.1 are satisfied.  Therefore, there exists a Lagrangian
isometric immersion $L_{f}:M^{n}\rightarrow {\bf H^{n}}$ whose
second fundamental form is given by
$h(X,Y)=I\sigma_{1}(X,Y)+J\sigma_{2}(X,Y)+K\sigma_{3}(X,Y)$. Up to
rigid motions of ${\bf H^{n}}$, it is unique by theorem 2.2.
However, if we put $e_{1}=\frac{1}{f}\frac{\partial}{\partial
x_{1}},e_{i}=\frac{\partial}{\partial x_{i}},i\geq 2$, then
$h(e_{1},e_{1})=\lambda_{1}Ie_{1}+\lambda_{2}Je_{1}+\lambda_{3}ke_{1},
\lambda_{i}=\frac{1}{f_{i}},i=1,2,3.$ and $h(e_{i},e_{j})=0$ for
$i,j\geq 2$.
\\[0.1in]
(b) The second fundamental form for $L$ is given by
\begin{eqnarray}\label{eq:2}
h(e_{1},e_{1})&=&\lambda_{1}Ie_{1}+\lambda_{2}Je_{1}+\lambda_{3}Ke_{1}\\\nonumber
h(e_{1},e_{j})&=&\mu_{1}Ie_{j}+\mu_{2}Je_{j}+\mu_{3}Ke_{j},j\geq2\\
\nonumber
h(e_{i},e_{i})&=&\mu_{1}Ie_{1}+\mu_{2}Je_{1}+\mu_{3}Ke_{1},i\geq2\\
\nonumber
h(e_{j},e_{k})&=&0,\hskip1pc j\neq k\geq 2, \nonumber
\end{eqnarray}
for some real valued functions $\lambda_{i}$,and $\mu_{i},i=1,2,3$
with respect to some orthonormal frame fields
$\{e_{1},e_{2},...,e_{n}\}$. Since $n\geq 3$, we can compute the
sectional curvature of the plane spanned by $e_{i}$, and $e_{j}$
for $i\neq j\geq 2$ which implies that
$\mu_{1}=\mu_{2}=\mu_{3}=0$. Now, the second fundamental form
given above becomes
\begin{eqnarray}\label{eq:3}
h(e_{1},e_{1})&=&\lambda_{1}Ie_{1}+\lambda_{2}Je_{1}+\lambda_{3}Ke_{1}\\\nonumber
h(e_{i},e_{j})&=&o\mbox{ for all }i,j\mbox{ except }i=1,\mbox{ and
}j=1\nonumber
\end{eqnarray}
By Codazzi equation, we get
\begin{eqnarray}\label{eq:4}
e_{i}(\lambda_{1})&=&\omega_{1}^{i}(e_{1})\lambda_{1}\\\nonumber
e_{i}(\lambda_{2})&=&\omega_{1}^{i}(e_{1})\lambda_{2}\\\nonumber
e_{i}(\lambda_{3})&=&\omega_{1}^{i}(e_{1})\lambda_{3} \nonumber
\end{eqnarray}
Also, we have
\begin{equation}\label{eq:5}
\nabla_{e_{i}}e_{1}=0.
\end{equation}
Let $\textit{D}$ and \textit{$D^{\bot}$} be the distributions
spanned by $e_{1}$ and $\{e_{2},...,e_{n}\}$, respectively. Since
\textit{D} is 1 dimensional, \textit{D} is integrable. Also,
\textit{$D^{\bot}$} is integrable because of (\ref{eq:5}).
Moreover, the leaves of $\textit{D}$ and \textit{$D^{\bot}$} are
totally geodesic submanifolds of ${\bf H^{n}}$. Since $\textit{D}$
and \textit{$D^{\bot}$} are integrable and perpendicular, there
exist local coordinates $\{x_{1},...,x_{n}\}$ such that
$\frac{\partial}{\partial x_{1}}$ spans \textit{D} and
$\{\frac{\partial}{\partial x_{2}},...,\frac{\partial}{\partial
x_{n}}\}$ spans \textit{$D^{\bot}$}. Since \textit{D} is 1
dimensional, we can choose $x_{1}$ such that
$\frac{\partial}{\partial
x_{1}}=|\lambda|e_{1},|\lambda|^{2}=\lambda_{1}^{2}+
\lambda_{2}^{2}+\lambda_{3}^{2}.$ Therefore, by Hiepko's
theorem(\cite{h}), $M$ is isometric to an open portion of the
twisted product manifold $_{f}I\times {\bf E^{n-1}}$ with the
twisted product metric
$g=f^{2}dx_{1}^{2}+\sum_{j=2}^{n}dx_{j}^{2},f=|\lambda|.$\\
We will consider the following case according to (\ref{eq:4}).\\
(\textbf{case} 1) $\omega_{i}^{i}(e_{1})=0$ for all $i\geq 2$\\
This condition implies that
$e_{j}(\lambda_{1})=e_{j}(\lambda_{2})=e_{j}(\lambda_{3})=0$ for
all $j\geq 2$. It means that the twistor function $f$ is a
function depending only on $x_{1}$. Using the twisted product
given above, we have
\begin{eqnarray} \label{eq:6}
\nabla_{\frac{\partial}{\partial x_{1}}}\frac{\partial}{\partial
x_{1}} & = & \frac{f'}{f}\frac{\partial}{\partial
x_{1}}\\\nonumber \nabla_{\frac{\partial}{\partial
x_{1}}}\frac{\partial}{\partial x_{i}} & = &
\nabla_{\frac{\partial}{\partial x_{1j}}}\frac{\partial}{\partial
x_{k}}=0,\mbox{ for }i,j,k=2,...,n \nonumber
\end{eqnarray}
Combining (\ref{eq:3}), (\ref{eq:6}) and Gauss' formula yield
\begin{eqnarray} \label{eq:7}
L_{x_{1}x_{1}} & = & (\frac{f'}{f}+(i\lambda_{1}+j\lambda_{2}+k\lambda_{3})|\lambda|)L_{x_{1}}\\
\nonumber L_{x_{1}x_{j}}& = & 0 \\
\nonumber L_{x_{j}x_{k}}& = & 0
\nonumber
\end{eqnarray}
By the third equation of (\ref{eq:7}), we get
$L(x_{1},...,x_{n})=D(x_{1})+\sum_{j=2}^{n}P_{j}(x_{1})x_{j}$\\
for some ${\bf H^{n}}$-valued functions $D,P_{2},...,P_{n}$. Using
the second equation of (\ref{eq:7}), we can find that
$P_{j}(x_{1})=c_{j}$ where $c_{j}$'s are constant vectors in ${\bf
H^{n}}$. Therefore, $L$ is a Lagrangian cylinder over a curve
$D=D(x_{1})$ whose rulings are $(n-1)$ plane parallel to
$x_{2}\cdots x_{n}$ plane in ${\bf H^{n}}$.\\
(\textbf{case} 2) $\omega_{i}^{i}(e_{1})\neq 0$ for some $i\geq 2$\\
If we assume $\lambda_{1},\lambda_{2},\lambda_{3}$ are positive,
then we have $e_{j}(\ln \lambda_{1})=e_{j}(\ln
\lambda_{2})=e_{j}(\ln \lambda_{3})$ for $j\geq 2$ and thus
$\frac{f}{\lambda_{1}},\frac{f}{\lambda_{2}},\frac{f}{\lambda_{3}}$
are all functions only depending on $x_{1}$. Using the twisted
product metric $g$, we get
\begin{eqnarray} \label{eq:8}
\nabla_{\frac{\partial}{\partial x_{1}}}\frac{\partial}{\partial
x_{1}} & = & \frac{f_{1}}{f}\frac{\partial}{\partial x_{1}}-
f\sum_{k=2}^{n}f_{k}\frac{\partial}{\partial x_{k}}\\
\nonumber
\nabla_{\frac{\partial}{\partial x_{1}}}\frac{\partial}{\partial
x_{i}} & = & \frac{f_{i}}{f}\frac{\partial}{\partial x_{1}},\\
\nonumber
\nabla_{\frac{\partial}{\partial
x_{j}}}\frac{\partial}{\partial x_{k}} & = & 0,\mbox{ for
}i,j,k=2,...,n \nonumber
\end{eqnarray}
By (\ref{eq:8}), we can compute the Riemannain curvature tensor
$R(\frac{\partial}{\partial x_{1}},\frac{\partial}{\partial
x_{j}})$ and using the fact that $M^{n}$ is flat, we have
$f_{jk}=0$ for all $j,k\geq 2$ which implies that $f$ is given by
$$f(x_{1},...,x_{n})=\beta(x_{1})+
\sum_{j=2}^{n}\alpha_{j}(x_{1})x_{j}$$ for some functions
$\beta,\alpha_{1},...,\alpha_{n}$ of $x_{1}$. \\
Therefore, by (a), $M^{n}$ is a Lagrangian submanifold described
in the statement (a).\\
More explicitly, we can describe the manifold $M^{n}$ as stated in
the theorem by doing the same computation done in \cite{c}.
\\[0.2in]
Next, we also have the following theorem for surfaces.\\
\begin{thm}\label{theorem:3.2}
Let $L:M^{2}\rightarrow {\bf H^{2}}$ be a Lagrangian $H$-umbilical
isometric immersion of a flat surface into ${\bf H^{2}}$ without
totally
geodesic points, then we have the following cases.\\
(a) $M^{2}$ is a Lagrangian cylinder in ${\bf H^{2}}$ i. e.
$L(x,y)=D(x)+cy$ for a curve $D=D(x)$ and a constant vector $c$ in
${\bf H^{2}}$.\\
(b) $L$ is given by $L(x,y)=D(x)+P(x)y$ for some ${\bf
H^{2}}$-valued curves $D$ and $P$. In fact, here $P=P(x)$ is a
special Legendre curve in $S^{7}\subset {\bf H^{2}}$.\\
(c) $L(x,y)=f(x)A(y)$, where $f$ is a ${\bf H}$-valued function
and $A$ is a curve in ${\bf H^{2}}$. In fact, it is a cone over a
curve $A$ in ${\bf H^{2}}$ plane.
\end{thm}
 \pf
We can have (a) and (b) if  $\mu_{1}=\mu_{2}=\mu_{3}=0$ in
(\ref{eq:2}). Its proof is same as theorem 3.1. For the case if
$|\mu|\neq 0$, we can also prove in the same way as theorem 3.1
and get the result (c).
\\[0.1in]
The author would like to thank Prof. B.Y. Chen for suggesting the
problem and for useful discussions on this topic.
%
%
%
%

%
%
%\address
Yun Myung Oh\\
Department of Mathematics\\
Michigan State University\\
E. Lansing, MI 48824\\
U.S.A.\\[0.1in]
Joon Hyuk Kang \\
Department of Mathematics\\
Andrews University\\
Berrien Springs, MI 49104\\
U.S.A.
\end{document}